\documentclass{amsart}
\usepackage{amssymb,latexsym,graphics,amsfonts}

\swapnumbers
 \theoremstyle{plain}

 \theoremstyle{definition}

 \numberwithin{equation}{subsection}

 \newcommand{\cal}[1]{\mathcal{#1}}

\newcommand{\cA}{\ensuremath{\mathcal A}}

\begin{document}

\title[Left introverted subspaces of duals of Banach algebras and...]
 {Left introverted subspaces of duals of Banach algebras and
 $WEAK^*-$continuous derivations on dual Banach algebras}

\author{ M. Eshaghi-Gordji  }

\address{Department of Mathematics, Semnan University, Semnan, Iran}
 \email{maj\_ess@yahoo.com and madjideg@walla.com}






\date{}

\dedicatory{}

\commby{}
\begin{abstract}
Let $X$ be a left introverted subspace of dual of a Banach algebra.
We study $Z_t(X^*),$ the topological center of Banach algebra $X^*$.
We fined the topological center of $(X\cA)^*$, when  $\cA$ has a
bounded right approximate identity and $\cA\subseteq X^*.$ So we
introduce a new notation of amenability for a dual Banach algebra
$\cal A$. A dual Banach algebra $\cal A$ is weakly Connes-amenable
if the first $weak^*-$continuous cohomology group of $\cal A$ with
coefficients in $\cal A$ is zero; i.e., $H^1_{w^*}(\cal A, \cal
A)=\{o\}$.  We study the  weak Connes-amenability of some dual
Banach algebras.

\end{abstract}
\maketitle
\section*{1. Introduction}
Let $\cA$ and be Banach algebra and let $X$ be a Banach
$\cA-$bimodule.  We can define right and left actions of $\cA$ on
the dual space $X^*$ of $X$ as follows
\begin{equation*}\quad \langle  fa, b \rangle = \langle f,
 ab \rangle ,
\quad \langle  af, b \rangle = \langle f,
 ba \rangle \quad (a,b \in \cal A, f\in X^* ).\end{equation*}
  Then $X ^{*}$ can be made into a
Banach $\cal A -$bimodule. For example, $\cal A$ itself is a Banach
$\cal A -$bimodule with respect to the product in $\cal A.$ Then
$\cal A^*$ is a Banach $\cal A -$bimodule.


The second dual space $\cal A ^{**}$ of a Banach algebra $\cal A $
admits a Banach algebra product known as first (left) Arens product.
We briefly recall the definition of this product.

For $m, n \in \cal A ^{**}$, their first (left) Arens product
indicated by $mn$ is given by
\begin{equation*}\langle m n, f \rangle = \langle m, n f \rangle \quad (f \in
\cal A ^{*}),\end{equation*}where $n f \in \cal A ^{*}$ is defined
by
\begin{equation*}\quad \langle n f, a \rangle = \langle n,
f a \rangle \quad (a \in \cal A )\hspace {0.5 cm}
[1].\end{equation*}

 Let $X$ be a Banach left $\cal
A-$module. Then we set $\cal A.X = \{ a.x : a \in \cal A , x \in X
\},$ and $\cal A X =\overline{lin \cal A.X}.$ If $\cA$ has a bounded
left approximate identity, then by Cohen-Hewitt factorization
theorem, (theorem 32.22 of [8]), $\cal A X =\overline{\cal A.X}.$ We
set similar definition for $X\cA$ when $X$ is a right $\cal A
-$module. A Banach algebra $\cal A$ is said to be dual if there is a
closed submodule ${\cal A}_*$ of $\cal A^*$ such that $\cal A={{\cal
A}_*}^*.$ Let $\cal A$ be a dual Banach algebra. A dual Banach $\cal
A$-bimodule $X$ is called normal if, for every $x\in X,$ the maps
$a\longmapsto a.x$ and $a\longmapsto x.a$ are $weak^*-$continuous
from $\cal A$ into $X$.  For example if $G$ is a locally compact
topological group, then $M(G)$ is a dual Banach algebra with predual
$C_0(G)$. Also if $\cal A$ is an Arens regular Banach algebra, then
$\cal A^{**}$ is a dual Banach algebra with predual $\cal A^*.$

 If $X$ is a Banach
$\cal A$-bimodule then a derivation from $\cal A$ into $X$ is a
linear map $D$, such that for every $a,b \in \cal A,$
$D(ab)=D(a).b+a.D(b).$ If $x\in X,$ and we define $\delta_x:\cal
A\longrightarrow X$ by $\delta_x(a)=a.x-x.a\hspace {0.5cm}(a\in \cal
A),$ then $\delta_x$ is a derivation. Derivations of this form are
called inner derivations. A Banach algebra $\cal A$ is amenable if
every bounded derivation from $\cal A$ into dual of every Banach
$\cal A$-bimodule $X$
 is inner; i.e., $H^1(\cal A, X^*)=\{o\}$ [9].
  Let $n\in \Bbb N$, then a Banach algebra $\cal A$
is n-weakly amenable if every (bounded) derivation from $\cal A$
into n-th dual of $\cal A$ is inner; i.e., $H^1(\cal A, {\cal
A}^{(n)})=\{o\}$ (see [4]). A dual Banach algebra $\cal A$ is
Connes-amenable if every $weak^*-$continuous derivation from $\cal
A$ into each normal dual Banach $\cal A$-bimodule $X$
 is inner; i.e., $H^1_{w^*}(\cal A, X)=\{o\}$, this definition was
introduced by V.Runde (see section 4 of [15]). In this paper we
study the $weak^*-$continuous derivations from $\cal A$ into itself
when $\cal A$ is a dual Banach algebra. We introduce the weak
Connes-amenability of dual Banach algebras as follows.
\paragraph{\bf Definition 1.1.} Let ${\cal A}$ be a dual Banach algebra.
Then $\cal A$ is said to be weakly Connes-amenable if every
$weak^*-$continuous derivation from $\cal A$ into $\cal A$ is inner;
i.e., $H^1_{w^*}(\cal A, \cal A)=\{o\}$.

Trivially we see that every Connes-amenable dual Banach algebra is
weakly Connes-amenable. We have already seen weakly Connes-amenable
dual Banach algebras which are not Connes-amenable.

\paragraph{\bf Example 1.2.} Let $\cal B$ be a Von-Neumann algebra. Then
$H^1_{w^*}(\cal B,\cal B)\subseteq H^1(\cal B,\cal B)=\{0\}$
(theorem 4.1.8 of [16]).  Thus $\cal B$ is weakly Connes-amenable
dual Banach algebra.
\paragraph{\bf Example 1.3.} Let $\cal A$ be a commutative semisimple
dual Banach algebra, then by commutative Singer-Warmer theorem, (see
for example [3; 18,16]) we have $H^1(\cal A, {\cal A})=\{o\}$, so
$\cal A$ is weakly Connes-amenable.

Let now $\cal A$ be a commutative Banach algebra which is Arens
regular and let $\cal A^{**}$ be semisimple. Trivially $\cal A^{**}$
is commutative. Then $\cal A^{**}$ is weakly Connes-amenable dual
Banach algebra.

Let $\cal A$ be a Banach algebra. The Banach $\cal A-$ submodule $X$
of $\cA^*$ is called left introverted if $\cal A^{**}X\subseteq X$
(i.e.$ X^*X\subseteq X$). Let $X$ be a left introverted Banach $\cal
A-$ submodule of $\cA^*$, then $X^*$ by the following product is a
Banach algebra:
$$\langle   x'y',x\rangle=\langle   x',y'\cdot x\rangle  \quad (x',y'\in X^*, x\in X) \quad
[1].$$ For each $ y' \in X^*$, the mapping $x' \longmapsto x' y'$ is
$weak^{*}-$continuous. However for certain $x'$, the mapping $y'
\longmapsto x' y'$ may fail to be $weak^{*}-$continuous. Due to this
lack of symmetry the topological center $Z_t(X^*)$ of $X^*$ is
defined by $$Z_t(X^*) := \{ x' \in X^*: y' \longmapsto x' \ y' :X^*
\longrightarrow X^* \hbox { ~~~is $weak^{*}-$continuous}\}
\hspace{0.5 cm} [12], [5].$$ If $X=\cA^*,$ then $Z_t(X^*)=Z_1(\cal A
^{**})$ is the left topological center of $\cA^{**}.$

 Let $X$ be a left introverted Banach $\cal A-$
submodule of $\cA^*$ and let $Y(\subseteq X)$ be a norm-closed
subspace of $X$. We denote $$Y^\perp:=\{m\in X^*: \langle m,f\rangle
=0 \; \hbox{~~for~~ every~}\; f\in Y\}.$$ Let $X$ be a Banach
$\cA-$bimodule. We set $Hom^r_{\cA}(X,X)$ the set of all elements
$T\in BL(X)$ such that $T(xa)=T(x)a, \hspace {0.5 cm}(x\in X, a\in
\cA).$ Similarly we set $Hom^l_{\cA}(X,X)$ the set of continuous
homomorphisms of left $\cA-$module $X$ to itself. If $X=\cA$, then
$Hom^l_{\cA}(X,X)$ is $LM(\cA)$ the set of all left multipliers of
$\cA.$  We denote $w^*-Hom^r_{\cA}(X^*,X^*)$ the set of all elements
in $Hom^r_{\cA}(X^*,X^*)$ which are $weak^*-$continuous. The paper
organized as follows. In section 2, we study the topological centers
of duals of left introverted subspaces of $\cA^*,$ when $\cA$ is a
Banach algebra with a bounded approximate identity. Indeed we fined
the topological center of $(X\cA)^*$, when   $\cA\subseteq X^*.$ We
study the relationship of the topological center of Banach algebra
$(\cA^*\cA)^*$ with the algebra $LM(\cA)$ of left multipliers of
$\cA.$ Of course we show that there is  an anti isomorphism between
$LM(\cA)$ and $Z_t((\cA^*\cA)^*)$ when $\cA Z_1(\cal A
^{**})\subseteq \cA$. We also provide an answer to the question when
the equality $wap(X)=X\cA$ is equivalent to the fact that $\cA$ is a
right ideal in $X^*$, when $X$ is a left introverted subspace of
$\cA^*.$ This includes results proved by Lau and \"{U}lger in [12].
 In
section 3, we study the weak Connes-amenability of dual Banach
algebras. So we give some examples of weakly Connes-amenable Banach
algebras which are not Connes-amenable.
\section*{2. Left introverted subspaces}
Let $X$ be a left introverted Banach $\cal A-$ submodule of $\cA^*$.
Then $X\cA$ is a left introverted subspace of $\cal A^*.$ We fined
the relations between the topological centers of $X^*$ and
$(X\cA)^*.$

 First we have the following results about left introverted subspaces.
\paragraph{\bf Lemma 2.1.}
Let $\cA$ be a Banach algebra and let  $X$ be a left introverted
subspace of $A^*$. Then the following assertions are equivalent.

(i) $X^*$ is a dual Banach algebra.

(ii) $Z_t(X^{*})=X^*$.

(iii) $\widehat X$ is a right $X^*$- submodule of $X^{**}$.

introverted subspace of $\cal B^*.$
\paragraph{\bf Proof.} $(i)\Longleftrightarrow (ii).$ It follows from 4.4.1 of [15].

$(ii)\Longrightarrow (iii).$ Let $x\in X$,  $x'\in X^*$  and let
$y'_\alpha~\stackrel{weak^*}{-\hspace{-.2cm}-\hspace{-.2cm}\longrightarrow}
y'$ in $X^*$. Then by (ii),
$x'y'_\alpha~\stackrel{weak^*}{-\hspace{-.2cm}-\hspace{-.2cm}\longrightarrow}
x'y'$ in $X^*$. So we have $$\langle \widehat{x}x',y'_\alpha
\rangle=\langle \widehat{x},x'y'_\alpha \rangle=\langle x'y'_\alpha,
x \rangle\longrightarrow \langle x'y',x \rangle=\langle
\widehat{x},x'y' \rangle=\langle \widehat{x}x',y' \rangle.$$ This
means that $\widehat{x}x':X^*\longrightarrow \Bbb C$ is
$weak^*-$continuous. Thus $\widehat{x}x'\in \widehat{X}.$

$(iii)\Longrightarrow (ii).$ Let $x'\in X^*$  and let
$y'_\alpha~\stackrel{weak^*}{-\hspace{-.2cm}-\hspace{-.2cm}\longrightarrow}
y'$ in $X^*$. Then for every $x\in X,$ we have
$$\langle x'y'_\alpha,
x \rangle=\langle y'_\alpha, xx' \rangle\longrightarrow \langle y',
xx' \rangle=\langle x'y',x \rangle.$$ Then (ii)
holds.\hfill$\blacksquare~$
\paragraph{\large\bf  Theorem 2.2.} Let $\cal A$ be a Banach algebra and let $Y\subseteq X$ be two
 norm-closed, $\cal A-$submodules of $\cal A^*$. If $X$ is left
 introverted, then

(i) The space $Y^\perp$ is a closed left ideal in $X^*$.

(ii) Suppose, further that $Y$ is left introverted. Then $Y^\perp$
is a closed ideal in $X^*.$

(iii) Suppose, further, that $X\cA\subseteq Y.$ Then $Y^\perp$ is a
left annihilator ideal in $X^*,$ and $Y^\perp\subseteq rad (X^*).$

(iv) $(X\cA)^*$ is, as a Banach algebra, isomorphic with the
quotient ${X^*}/{(X\cA)^\perp}.$

(v) Suppose that $X \cA\neq X$. Then $rad(X^*)\neq \{0\}.$
\paragraph{\bf  Proof.} Let $y'\in Y^\perp$ and $x'\in X^*.$ Then

(i) For each $y\in Y,$ we have $\langle x'y',y\rangle=\langle x',y'y
\rangle=\langle x',0 \rangle =0$ because $Y$ is  $\cal A-$submodules
of $\cal A^*$, and for each $a\in \cA$, we have $\langle
y'y,a\rangle=\langle y',ya \rangle=0.$

(ii) For each  $y\in Y,$ we have $x'y \in Y.$ Then $\langle
y'x',y\rangle=\langle y',x'y \rangle=0.$

(iii) For $x\in X$ and $a\in \cA,$ we have $\langle
y'x,a\rangle=\langle y',xa \rangle=0.$ Thus $X'Y'=0$ in $X^*$. Then
$Y^\perp$ is a left annihilator ideal in $X^*$. So $Y^\perp
\subseteq X^*.$

(iv), (v) follow from (iii).\hfill$\blacksquare~$

Let $X$ be a left introverted Banach $\cal A-$ submodule of $\cA^*$,
$x'$ be an element in $X^*$ and $m\in (X\cA)^*.$ For $x\in X$ we
define $\widehat {x.x'}:X^*\longrightarrow \Bbb C$ by $\langle
\widehat {x.x'},y'\rangle=\langle x'y',x\rangle$ and for each $z\in
X\cA$, we define $\widetilde{z.m}:(X\cA)^*\longrightarrow \Bbb C$ by
$\langle \widetilde{z.m},y'\rangle=\langle my',z\rangle.$ For $x'\in
X^*$ we take $\check{x'}=x'\mid_{X\cA}.$  We have
$am.n=am.\check{n},$ for each $a\in \cA.$
\paragraph{\large\bf  Theorem 2.3.} Let $\cal A$ be a Banach algebra
with a bounded right approximate identity,
 and let $X$ be a norm-closed, left introverted  $\cal A-$submodules of $\cal
 A^*$. If $\cA \subseteq X^*,$ then

(i) $(X\cA)^\perp$ is the ideal of right annihilators in $X^*$.

(ii) $(X\cA)^*$ is isomorphic (as a Banach algebra) with
$Hom^r_{\cA}(X,X)$.

(iii) $X^*\cong (X\cA)^* \oplus(X\cA)^\perp$ as Banach spaces.

(iv) Let $m\in X^*,$ then $m\in Z_t((X\cA)^*)$ if and only if for
each $a\in \cA,$ $am\in Z_t(X^*).$

(v) If $X=Y^*$ for Banach $\cA-$bimodule $Y$ and
$\cA(Z_t(X\cA)^*)\subseteq \cA$, then there is an isomorphism
between $Z_t((X\cA)^*)$ and $w^*-Hom^r_{\cA}(X,X)$.

(vi) By conditions in (v), there is an anti isomorphism between
$Z_t((X\cA)^*)$ and $Hom^l_{\cA}(Y,Y)$.
\paragraph{\bf  Proof.} (i) Since $X$ is left introverted, then $X\cA$ is a left
introverted subspace of $\cA^*.$ By above theorem, $(X\cA)^\perp$ is
contained in the left annihilators of $X^*.$ On the other hand, if
$x'\in X^*$ is a right annihilator then $ax'=0$ for each $a\in \cA$
because $\cA \subseteq X^*.$ Then $\langle x',xa\rangle=\langle
ax',x\rangle=0$ for all $x\in X, a\in \cA,$ so that $x'\in
(X\cA)^\perp.$

(ii) We define $\phi:X^*\longrightarrow Hom^r_{\cA}(X,X),$ by
$\phi(x')x=x'x$ for all $x'\in X^*, x\in X$. Obviously for each
$x'\in X^*$, $\phi(x')$ is a right $\cA-$module homomorphism of $X$.
 The kernel of $\phi$ is $(X\cA)^\perp$. Since $\cal A$ has bounded right approximate
 identity, then $\cA^{**}$ (and then $X^*$) has a right identity.
 Let  $E$ be a right identity of $X^*$. For each $T\in
 Hom^r_{\cA}(X,X),$ we define $T_{\phi}\in X^*$ by $\langle T_{\phi},x\rangle=\langle E,Tx \rangle.$
For each $x\in X, a\in \cA,$ we have \begin{eqnarray*}\langle
T_{\phi}x,a\rangle=\langle T_{\phi},xa\rangle=\langle E,T(xa)\rangle
=\langle E,(Tx)a\rangle=\langle aE,(Tx)\rangle= \langle
a,Tx\rangle=\langle Tx,a\rangle. \end{eqnarray*}
 Thus
$Tx=\phi(T_{\phi})x.$ This shows that $\phi$ is onto, so by above
theorem, $Hom^r_{\cA}(X,X)\cong {X^*}/{(X\cA)^\perp}\cong (X\cA)^*$
.

(iii) Let $E$ be as above, then for each $T\in
 Hom^r_{\cA}(X,X),$ we have $\langle ET_{\phi},x\rangle=
 \langle E,T_{\phi}x \rangle=\langle E,Tx \rangle=\langle
 T_{\phi},x\rangle.$ Then $ET_{\phi}=T_{\phi}.$ Thus $EX^*\cong
 (X\cA)^*,$ and $(X\cA)^\perp=\{x'-Ex':x'\in X^*\}.$ So we have
 $X^*\cong (X\cA)^* \oplus(X\cA)^\perp.$

(iv) Let $m\in Z_t((X\cA)^*)$ and let $a\in \cA$, then for each
$x\in X,$ $\widetilde{xa.m}:(X\cA)^*\longrightarrow \Bbb C$  is
$weak^*-$continuous. Then $\widetilde{xa.m}\in X\cA$. For $n\in X^*$
we have $am.n=am.\check{n},$ then $$\langle
\widehat{x.am},n\rangle=\langle x,am.n\rangle=\langle
x,am.\check{n}\rangle =\langle \widetilde{xa.m},\check{n}\rangle
=\langle \widetilde{xa.m},n\rangle \hspace {0.5 cm} (2.1).$$ Thus
$\widehat{x.am}=\widetilde{xa.m}\subseteq X\cA\subseteq X.$ Then
$\widehat{x.am}:X^*\longrightarrow \Bbb C$ is $weak^*-$continuous.
It means that $am\in Z_t(X^*).$ Conversely let $m\in X^*,$ satisfies
for each $a\in \cA,$ $am\in Z_t(X^*).$  We show that for each $x\in
X, a\in \cA$, $\widetilde{xa.m}:(X\cA)^*\longrightarrow \Bbb C$  is
$weak^*-$continuous;(i.e. $m\in Z_t((X\cA)^*)$). Since $am\in
Z_t(X^*),$ then $\widehat{x.am}:X^*\longrightarrow \Bbb C$ is
$weak^*-$continuous. By (2.1), we have
$\widehat{x.am}=\widetilde{xa.m}$ and the proof is complete.

(v) By (ii) we know that $T\longmapsto
T_{\phi}:Hom^r_{\cA}(X,X)\longrightarrow (X\cA)^* $ is an
isomorphism. We denote this map with $\Psi.$ Let $a\in \cA$ then for
each $x\in X$ we have
\begin{eqnarray*}\langle a\Psi(T),x\rangle=\langle
aT_{\phi},x\rangle=\langle T_{\phi},xa\rangle =\langle
E,T(xa)\rangle=\langle E,(Tx)a\rangle =\langle Tx,a\rangle \hspace
{0.5 cm}(2.2). \end{eqnarray*}
 Now, let $T\in w^*-Hom^r_{\cA}(X,X)$ and let
$y'_\alpha~\stackrel{weak^*}{-\hspace{-.2cm}-\hspace{-.2cm}\longrightarrow}
y'$ in $X$, then
$Ty'_\alpha~\stackrel{weak^*}{-\hspace{-.2cm}-\hspace{-.2cm}\longrightarrow}
Ty'$ in $X$, and by (2.2), we have $$\lim _{\alpha}\langle
a\Psi(T),y'_{\alpha}\rangle=\lim _{\alpha}\langle
Ty'_{\alpha},a\rangle=\langle Ty',a\rangle=\langle
a\Psi(T),y'\rangle.$$ This means that $a\Psi(T):X=Y^*\longrightarrow
\Bbb C$ is $weak^*-$continuous. Then $a\Psi(T)\in Y$. Let now
$x'_\alpha~\stackrel{weak^*}{-\hspace{-.2cm}-\hspace{-.2cm}\longrightarrow}
x'$ in $X^*$, then for each $x\in X$ we have \begin{eqnarray*} \lim
_{\alpha}\langle a\Psi(T)x'_{\alpha},x\rangle=\lim _{\alpha}\langle
a\Psi(T),x'_{\alpha}x\rangle = \lim _{\alpha}\langle
x'_{\alpha}x,a\Psi(T)\rangle =\langle x'x,a\Psi(T)\rangle=\langle
a\Psi(T)x',x\rangle. \end{eqnarray*}
 Thus
$a\Psi(T)x'_\alpha~\stackrel{weak^*}{-\hspace{-.2cm}-\hspace{-.2cm}\longrightarrow}
a\Psi(T)x'$ in $X^*,$ so $a\Psi(T)\in Z_t(X^*).$ Therefore by (iv)
above, $\Psi(T) \in Z_t((X\cA)^*).$ Conversely let  $T\in
 Hom^r_{\cA}(X,X)$ and $\Psi(T) \in Z_t((X\cA)^*).$ We will show
 that $T\in w^*-Hom^r_{\cA}(X,X)$. Let $x_\alpha~\stackrel{weak^*}{-\hspace{-.2cm}-\hspace{-.2cm}\longrightarrow}
x$ in $X,$ and let $a\in \cA,$ then $a\Psi(T)\in \cA$. Thus
\begin{eqnarray*}\lim _{\alpha}\langle Tx_{\alpha},a\rangle&=&\lim _{\alpha}\langle
\Psi(T),x_{\alpha}a\rangle=\lim _{\alpha}\langle
a\Psi(T),x_{\alpha}\rangle= \lim _{\alpha}\langle
x_{\alpha},a\Psi(T)\rangle\\
&=& \langle x,a\Psi(T)\rangle= \langle \Psi(T),xa\rangle=\langle
Tx,a\rangle.
\end{eqnarray*}
 Therefore $T$ is
$weak^*-$continuous.

(vi) Since $X=Y^*$, then $T\in
 BL(X,X)$ is $weak^*-$continuous if and only if
 there exists $S\in
 BL(Y,Y)$ such that $T=S^*$ the adjoint of $S.$ Let now $T\in
 Hom^r_{\cA}(X,X), a\in \cA, y'\in Y^*=X.$ Then $T(y'a)=T(y')a$ so
 for each $y\in Y,$ we have
 \begin{eqnarray*}\langle y',aS(y)\rangle=\langle y'a,S(y)\rangle=\langle T(y'a),y\rangle
=
 \langle T(y')a,y\rangle=\langle T(y'),ay\rangle
 =\langle y',S(ay)\rangle.
 \end{eqnarray*} Thus by Hahn-Banach theorem, we have
 $S(ay)=aS(y)$, and $S\in Hom^l_{\cA}(Y,Y).$ Now by (v), it is easy to show
 that the map $\Psi o~~~*:Hom^l_{\cA}(Y,Y)\longrightarrow
 Z_t((X\cA)^*)$ is an anti isomorphism.\hfill$\blacksquare~$

\paragraph{\bf Corollary 2.4.} (corollary 1.2 of [2]) Let $\cA$ be a
Banach algebra with a bounded right approximate identity. Then we
have $\cA^{**}\cong (\cA^*\cA)^* \oplus(\cA^*\cA)^\perp$ as Banach
spaces.
\paragraph{\bf Corollary 2.5.} Let $\cA$ be a Banach algebra with a
bounded right approximate identity, and let $X$ be a left
introverted $\cA-$submodule of $\cA^*.$ Then $\cA Z_t((X\cA)^*)=\cA
Z_t((X)^*).$
\paragraph{\bf Proof.} Since $\cA$ has a bounded right approximate
identity, then $\cA \cA=\cA$. By (iv) of the preceding theorem, $\cA
Z_t((X\cA)^*)\subseteq  Z_t((X)^*).$ Thus $\cA Z_t((X\cA)^*)=\cA \cA
Z_t((X\cA)^*)\subseteq \cA Z_t((X)^*).$ Now let $x'\in Z_t((X)^*)$
and let $n$ be the restriction of $x'$
 to $X\cA$. Then it is easy to show that $n\in Z_t((X\cA)^*)$
and that $ax'=an.$ Hence $\cA Z_t((X\cA)^*)=\cA
Z_t((X)^*).$\hfill$\blacksquare~$
\paragraph{\bf Corollary 2.6.} Let $\cA$ be a Banach algebra with a
bounded right approximate identity. If $\cA Z_1(\cal A
^{**})\subseteq \cA$ , then there exists an anti isomorphism between
$LM(\cA)$ and $Z_t((\cA^*\cA)^*)$.

Let $\cA$ be a Banach algebra, and let $X$ be a left introverted
$\cA-$submodule of $\cA^*.$ We denote by $wap(X)$ the set of
elements $x$ in $X$ for which  $\widehat {x.x'}:X^*\longrightarrow
\Bbb C$ is $weak^*-$continuous for each  $x'\in X^*$.  Obviously we
have the following assertions.

(i) $wap(X)$ is an $\cA$-submodule of $X$.

(ii) $wap(X)=X$ if and only if $Z_t((X)^*)=X^*.$
\paragraph{\bf Theorem 2.7.}
Let $\cA$ be a Banach algebra with a bounded approximate identity,
and let $X$ be a left introverted $\cA-$submodule of $\cA^*.$ If
$\cA \subseteq X^*,$ then $wap(X)$ is essential $\cA$-bimodule, and
$$wap(X)\subseteq X.\cA \cap \cA.X. $$
\paragraph{\bf Proof.} Let $(e_\alpha)_\alpha$ be a bounded approximate
identity for $\cA$ with bound M. Let $x\in wap(X)$, it is easy to
show that the set $\{x.a: a\in \cA, \|a\|\leq M \}$ is relatively
weakly compact in $X$. Let $U=\overline{\{a.x:\|a\|\leq M , a\in
\cA\}}^w.$ Then $U$ is compact. So we may suppose that the net
$(e_\alpha .x)_\alpha$ converges in the weak topology of $X$ in $U$.
Let $e_\alpha
.x~\stackrel{weakly}{-\hspace{-.2cm}-\hspace{-.2cm}\longrightarrow}
u$ in $U.$ Since $\cA \subseteq X^*,$ then for each $a\in \cA$, we
have $$\langle u,a \rangle=\lim_\alpha\langle \widehat{a},e_\alpha
.x\rangle=\lim_\alpha\langle x,
 ae_\alpha\rangle=\langle x,a \rangle.$$
Therefore $x=u$. It follows that $x\in \cA.X.$ Similarly,  $x\in X.
\cA .$ Thus $wap(X)$ is essential $\cA$-bimodule. By Cohen$^,$s
factorization theorem, the result follows.\hfill$\blacksquare~$
\paragraph{\bf Theorem 2.8.}
Let $\cA$ be a Banach algebra with a bounded approximate identity,
and let $X$ be a left introverted $\cA-$submodule of $\cA^*.$ Then
the following assertions are equivalent.

(i) $X\cA\subseteq wap(X).$

(ii) $\cA X^*\subseteq Z_t(X^*).$

(iii) $\cA X^*\subseteq \cA Z_t((X\cA)^*).$

(iv) $Z_t((X\cA)^*)=(X\cA)^*.$

(v) $(X\cA)^*$ is a dual Banach algebra.
\paragraph{\bf Proof.} $(i)\Longrightarrow (ii).$ Let $a\in \cA, x'\in
X^*,$ and let
$y'_\alpha~\stackrel{weak^*}{-\hspace{-.2cm}-\hspace{-.2cm}\longrightarrow}
y'$ in $X^*$. We have to show that
$a.x'y'_\alpha~\stackrel{weak^*}{-\hspace{-.2cm}-\hspace{-.2cm}\longrightarrow}
a.x'y'$ in $X^*$. To this end, let $x\in X,$ since $xa\in wap(X),$
then we have
$$\lim_\alpha \langle a.x'y'_\alpha,
x \rangle=\lim_\alpha \langle \widehat {xax'},y'_\alpha
\rangle=\langle \widehat{xax'},y' \rangle=\langle a.x'y',x
\rangle.$$ So the result follows. \\
$(ii)\Longrightarrow(iii).$ Let $\cA X^*\subseteq Z_t(X^*).$ Since
$\cA$ has a bounded approximate identity, then by corollary 5.2 and
Cohen$^,$s factorization theorem, we have $\cA X^*=\cA\cA
X^*\subseteq \cA Z_t(X^*)= \cA Z_t((X\cA)^*).$

 $(iii)\Longrightarrow
(iv).$ Let $m\in (X\cA)^*$ and let
$n_\alpha~\stackrel{weak^*}{-\hspace{-.2cm}-\hspace{-.2cm}\longrightarrow}
n$ in $(X\cA)^*$. We have to show that
$mn_\alpha~\stackrel{weak^*}{-\hspace{-.2cm}-\hspace{-.2cm}\longrightarrow}
mn$ in $(X\cA)^*$. To this end let $a\in \cA, x\in X,$ and let
$\tilde{m}$ be the Hahn-Banach extension of $m$ to $X$. Then
$a\tilde{m}\in \cA X^*\subseteq \cA Z_t((X\cA)^*).$ Thus there are
$b\in \cA, z\in \cA Z_t((X\cA)^*)$ such that $a\tilde{m}=bz$.
Therefore we have
\begin{eqnarray*}\lim _{\alpha}\langle mn_{\alpha},xa\rangle &=&\lim _{\alpha}\langle
amn_{\alpha},x\rangle=\lim _{\alpha}\langle
a\tilde{m}n_{\alpha},x\rangle=\lim _{\alpha}\langle
bzn_{\alpha},x\rangle= \lim _{\alpha}\langle zn_{\alpha},xb\rangle
\\
&=&\langle zn,xb\rangle=\langle bzn,x\rangle=\langle
a\tilde{m}n,x\rangle=\langle mn,xa\rangle.
\end{eqnarray*}
$(iv)\Longrightarrow (i).$ Let $a\in \cA, x\in X,$ and  $x'\in X^*.$
If
$y'_\alpha~\stackrel{weak^*}{-\hspace{-.2cm}-\hspace{-.2cm}\longrightarrow}
y'$ in $X^*$. Then
$\check{y'}_\alpha~\stackrel{weak^*}{-\hspace{-.2cm}-\hspace{-.2cm}\longrightarrow}
\check{y'}$ in $(X \cA)^*$. since $\check{x'} \in (X\cA)^*=
Z_t((X\cA)^*),$ then we have
$$\lim_\alpha \langle \widehat{xax'},y'_\alpha \rangle=\lim_\alpha \langle
x'y'_\alpha, xa\rangle=\lim_\alpha \langle
\check{x'}\check{y'}_\alpha, xa\rangle=\langle \check{x'}\check{y'},
xa\rangle=\langle x'y', xa\rangle=\langle
\widehat{xax'},y'\rangle.$$ Thus $xa\in wap(X).$\\
$(iv)\Longleftrightarrow (v).$ It follows from lemma
2.1.\hfill$\blacksquare~$
\paragraph{\bf Corollary 2.9.} Let $\cA$ be a Banach algebra with a
bounded approximate identity, and let $X$ be a left introverted
$\cA-$submodule of $\cA^*.$ If $\cA \subseteq X^*,$ then the
following assertions are equivalent.

(i) $\cA $ is a right ideal of $X^*$.

(i) $\cA $ is a right ideal of $(X \cA)^*$.

(iii) $wap(X)=X\cA$ and $\cA Z_t(X^*)\subseteq \cA.$
\section*{3. weak Connes-amenability}
In this section we study the first $weak^*-$continuous cohomology
group of $\cal A$ with coefficients in $\cal A$, when $\cA$ is a
dual Banach algebra. Indeed we show that an Arens regular Banach
algebra $\cal A$ is 2-weakly amenable if and only if the second dual
of $\cal A$ is weakly Connes-amenable. So we prove that a dual
Banach algebra $\cal A$ is weakly Connes-amenable if it is 2-weakly
amenable.

 Let $G$ be a locally
compact topological, inner amenable group, then the dual Banach
algebra $M(G)$ is Connes-amenable if and only if $L^1(G)$ is
amenable (see section 4 of [15]). Also $L^1(G)$ is always weakly
amenable (see [10] or [6]). In the following we show that $M(G)$ is
always weakly Connes-amenable.
\paragraph{\bf Theorem 3.1.} Let $X$ be a left introverted $\cal A-$submodule of $L^\infty(G)$ such that
$C_0(G)\subseteq X\subseteq CB(G)$. Let $Z_t(X^{*})=X^*$, and let
$L^1(G)$ be a closed ideal of $X^*$. Then $X^*$ is weakly
Connes-amenable.
\paragraph{\bf Proof.} We know that for every (bounded) derivation
$D:L^1(G)\longrightarrow L^1(G),$ there is $\mu \in M(G)$ such that
for every $a\in L^1(G),$ $D(a)=a\mu -\mu a$ [13; corollary 1.2].
 Let  $D:X^*\longrightarrow X^*$ be a $weak^*-$continuous derivation,
since $L^1(G)$ is a two sided ideal in $X^*$, then for every $a,b
\in L^1(G),$ we have $D(ab)=D(a).b+a.D(b)\in L^1(G).$ The proof of
the converse is easy because $L^1(G)$ is $weak^*-$dense in $X^*.$
\hfill$\blacksquare~$
\paragraph{\bf Corollary 3.2.} For every locally compact topological group
$G$, $M(G)$ is weakly Connes-amenable.
\paragraph{\bf Theorem 3.3.} Let $\cA$ be a commutative, semisimple  Banach algebra with a bounded approximate identity,
 and let  $X$ be a left introverted
subspace of $A^*$. If $ Z_t(X^{*})=X^*$ and $X^*$ contains $\cal A$
as a two sided ideal, then $X^*$ is weakly Connes-amenable.
\paragraph{\bf Proof.}
 Let  $D:X^*\longrightarrow X^*$ be a $weak^*-$continuous derivation,
 we can show that  $D(\cal A)\subseteq \cal A$. Therefore by
 commutative Singer-Warmer theorem [3; 18,16],
there is a $a \in \cal A(\subseteq X^*)$ such that  $D\mid_{\cal
A}=\delta_a.$ Since $\cal A$ is $weak^*-$dense in $X^*,$ and $D$ is
$weak^*-$continuous, then $D=\delta _a.$\hfill$\blacksquare~$

Let $\cA$ be a commutative, semisimple  Banach algebra with a
bounded approximate identity. If $\cA$ is an ideal in $\cA^{**},$
then by above theorem, the algebra $(\cA^*\cA)^*$ is weakly
Connes-amenable.
\paragraph{\bf Theorem 3.4.} Let $\cA$ be a Banach algebra and let  $X$ be a left introverted
$\cal A-$submodule of $A^*$. If $\widehat{\cal A}\subseteq
Z_t(X^{*})=X^*$, then the following assertions are equivalent.

(i) $X^{*}$ is weakly Connes-amenable.

(ii) $H^1(\cal A, X^*)=\{0\}.$
\paragraph{\bf Proof.} $(i)\Longrightarrow (ii)$. Let $D:{\cal
A}\longrightarrow X^{*}$ be a (bounded) derivation. Then by
proposition 1.7 of [4], we know that $D^{**}:{\cal
A}^{**}\longrightarrow ({X}^{*})^{**}$ the second transpose of $D$
is a derivation. We define $D_1:{X}^{*}\longrightarrow {X}^{*}$ by

$\hspace {3cm}\langle D_1(x'),x\rangle=\langle
D^{**}(x'),\widehat{x}\rangle \hspace {1cm}
 (x'\in {X}^{*}, x \in X).$\\
  Since $Z_t(X^{*})=X^*$, then by lemma 2.1, $\widehat X$ is a  $X^*$- submodule of $X^{**}$. Then for every
$x',y'\in X^{*}$ and  $x\in \cal A^{*}$, we have
\begin{align*}\langle D_1(x'y'),x \rangle&=\langle
D^{**}(x'y'),\widehat{x}\rangle=\langle
D^{**}(x')y',\widehat{x}\rangle+\langle
x'D^{**}(y'),\widehat{x}\rangle\\
&=\langle D^{**}(x'),y'\widehat{x}\rangle+\langle
D^{**}(y'),\widehat{x}x'\rangle=\langle
D^{**}(x'),\widehat{y'{x}}\rangle+\langle
D^{**}(y'),\widehat{xx'}\rangle\\
&=\langle D_1(x'),y'{x}\rangle+\langle D_1(y'),xx'\rangle= \langle
D_1(x')y',{x}\rangle+\langle x'D_1(y'),x\rangle.
\end{align*}
Then $D_1$ is a derivation.  Now let
$x'_\alpha~\stackrel{weak^*}{-\hspace{-.2cm}-\hspace{-.2cm}\longrightarrow}
x'$ in $X^*$. Since $D^{**}$ is $weak^*-$continuous, then for every
$x\in X$, we have
$$\lim _\alpha \langle D_1(x'_\alpha) ,x\rangle=\lim _\alpha \langle
D^{**}(x'_\alpha) ,\widehat{x}\rangle = \langle D^{**}(x')
,\widehat{x}\rangle=\langle D_1(x') ,x\rangle.$$ This means that
$D_1$ is $weak^*-weak^*-$continuous. Then there exists $x'\in X^*$
such that $D_1=\delta _{x'},$ so $D=\delta _{x'}$.

$(ii)\Longrightarrow (i)$. Let $D:X^{*}\longrightarrow X^{*}$ be a
$weak^*-$continuous derivation, then  $D\mid _{\cal A}:{\cal A}
\longrightarrow X^{*}$ is a bounded derivation. Thus there is $x'\in
X^{*}$ which $D(\widehat{a})=\widehat{a}x'-x'\widehat {a}$ for every
$a\in \cal A.$ Since $X^*$ is a dual Banach algebra, then $\delta
_{x'}:X^{*}\longrightarrow X^{*}$ is $weak^*-$continuous. On the
other hand  $\widehat {\cal A}$ is $weak^*-$dense in $X^{*},$ and
$D$ is $weak^*-$continuous, then we have $D=\delta
_{x'}.$\hfill$\blacksquare~$
\paragraph{\bf Corollary 3.5.}
Let $\cal A$ be an Arens regular Banach algebra, then  $\cal A^{**}$
is weakly Connes-amenable if and only if $\cal A$ is 2-weakly
amenable.
\paragraph{\bf Theorem 3.6.} Let $\cal A$ be a  dual  Banach
algebra, then the following assertions hold.

(i) If  $\cal A$ is 2-weakly amenable, then $\cal A$ is weakly
Connes-amenable.

(ii) If $\cal A$ is  Arens regular and  $\cal A^{**}$ is weakly
Connes-amenable, then $\cal A$ is weakly Connes-amenable.
\paragraph{\bf Proof.} (i) Let $\cal A$ be a dual algebra with predual
$\cal A_*$, and let $D:{\cal A}\longrightarrow {\cal A}$ be a
$weak^*-$continuous derivation, then $D$ is bounded. The natural
embedding $\widehat{} : \cal A\longrightarrow \cal A^{**}$ is an
$\cal A-$bimodule morphism, then $\widehat{}~~~ oD : \cal
A\longrightarrow \cal A^{**}$ is a bounded derivation. Since $\cal
A$ is 2-weakly amenable, then there exists $a''\in \cal A^{**}$ such
that $\widehat{}~~~ oD=\delta_{a''}$. We have the following direct
sum decomposition
$$\cal A^{**}=\cal A\oplus {\cal A_*}^\perp$$
as $\cal A-$bimodules [7]. Let $\pi :\cal A^{**}\longrightarrow \cal
A$ be the projection map. Then $\pi$ is an $\cal A-$bimodule
morphism, and we result that $D=\delta_{\pi(a'')}$. (ii) follows
from corollary 3.5 and (i). \hfill$\blacksquare~$

In the following (example 1) we will show that the converse of (i)
(or (ii)) in above theorem dos not holds.
\paragraph{\bf Examples.}

1- Let $\omega: \Bbb Z\longrightarrow \Bbb R$ define by
$\omega(n)=1+\left\vert n\right\vert$ and let
$$l^1(\Bbb Z,\omega)=\{\sum _{n\in \Bbb
Z}f(n)\delta_n:\left\Vert\,f\right\Vert_{\omega}=\sum \left\vert
f(n)\right\vert \omega(n)<\infty\}.$$ Then $l^1(\Bbb Z,\omega)$ is a
Banach algebra with respect to the convolution product defined by
the requirement that
$$\delta_m \delta_n=\delta_{mn} \hspace {1cm} (m,n \in \Bbb Z).$$
We define $$l^\infty(\Bbb Z, \frac{1}{\omega})=\{\lambda=\sum _{n\in
\Bbb Z}\lambda(n)\lambda_n:sup \frac{\left\vert
\lambda(n)\right\vert}{\omega(n)}<\infty\},$$ and
$$C_0(\Bbb Z,
\frac{1}{\omega})=\{\lambda \in l^\infty(\Bbb Z, \frac{1}{\omega})
:\frac{\left\vert \lambda\right\vert}{\omega(n)}\in C_0(\Bbb Z)\}.$$
Then $\cal A=l^1(\Bbb Z,\omega)$ is an Arens regular dual Banach
algebra with predual $C_0(\Bbb Z, \frac{1}{\omega})$ [5]. $\cal A$
is commutative and semisimple, then $\cal A$ is weakly
Connes-amenable (see example 1.3). On the other hand we know that
$\cal A$ is not 2-weakly amenable [5], (i.e. $\cal A^{**}$ is not
weakly Connes-amenable).


2- We propose that the Banach algebra of approximate operators $\cal
A (E)$ is defined as the closure of $\cal F (E)$ in $(\cal B
(E),\|.\|).$ Then by lemma 5.5 of [4], $\cal A (E)$ is Arens
regular. By applying corollary 3.5 above and corollary 5.10 of [4],
we have the following.

 (i) If $dim K_{E^*}\leq 1$ and
$E^*$ has the Radon-Nikodym property, then  $\cal A (E)^{**}$ is
weakly Connes-amenable dual Banach algebra.

 (ii) If $dim K_{E^*}\geq
1$, then $\cal A (E)^{**}$ is not weakly Connes-amenable.

 Let $E$ be a reflexive space with approximate
property, then $\cal N(E)$ the algebra of nuclear operators on $E$,
is Arens regular. Then $\cal N(E)^{**}$ is a dual Banach algebra. As
in corollary 5.4 of [4], we know that $\cal N(E)$ is not 2-weakly
amenable. Thus by corollary 3.5 above,  $\cal N(E)^{**}$ is not
weakly Connes-amenable.

3- The algebra $C^{(1)}(\Bbb I)$ consists of the continuously
differentiable functions on the unit interval $\Bbb I=[0,1]$;
$C^{(1)}(\Bbb I)$ is a Banach function algebra on $\Bbb I$ with
respect to the norm $\|f\|_1=\|f\|_{\Bbb I}+\|f'\|_{\Bbb I}$ $(f\in
C^{(1)}(\Bbb I)).$ By proposition 3.3 of [4], $C^{(1)}(\Bbb I)$ is
Arens regular but it is not 2-weakly amenable. Thus by corollary 3.5
above, $C^{(1)}(\Bbb I)^{**}$ is a dual Banach algebra which is not
weakly Connes-amenable.

4- For a function $f\in L^1(\Bbb T),$ the associated Fourier series
is $(\hat {f}(n):n\in \Bbb Z).$ For $\alpha \in (0,1)$ the
associated Beurling algebra $A_\alpha (\Bbb T)$ on $\Bbb T$ consists
of the continuous functions $f$ on $\Bbb T$ such that
$\|f\|_\alpha=\sum_{n\in \Bbb Z}\mid {\hat {f}(n)}\mid{(1+\mid
n\mid)}^{\alpha}<\infty.$ By proposition 3.7 of [4], ${A_\alpha
(\Bbb T)}$ is Arens regular and 2-weakly amenable. Then by applying
corollary 3.5 above, ${A_\alpha (\Bbb T)}^{**}$ is weakly
Connes-amenable.

\end{document}